\documentclass{amsart}%
\usepackage{amsmath,amssymb,amsthm}
\usepackage[latin1]{inputenc}
\usepackage{tabularx,multicol}
\usepackage{enumerate}
\usepackage{amsmath}
\usepackage{amsfonts}
\usepackage{amssymb}
\usepackage{graphicx}
\usepackage{latexsym}%
\theoremstyle{plain}
\newtheorem{theo}{Theorem}

\newtheorem{proposition}[theo]{Proposition}
\newtheorem{lemma}[theo]{Lemma}

\theoremstyle{definition}
\newtheorem{rem}[theo]{Remark}

\begin{document}
\title[Single point source SBM]{On the large scale behavior of super-Brownian motion\\in three dimensions with a single point source}
\date{\today\quad SinglePoint82.tex\quad WIAS Preprint No.\ 1154 of July 21,
2005\quad ISSN 0946\thinspace--\thinspace8633\quad Running head:\quad Single
point source SBM\quad Corresponding author: Klaus Fleischmann}
\author{Klaus Fleischmann}
\address{Weierstra\ss \ Institute for Applied Analysis and Stochastics, Mohrenstrasse
39, D-10117 Berlin, Germany}
\email{fleischm@wias-berlin.de}
\author{Carl Mueller}
\address{Department of Mathematics, University of Rochester, Rochester, NY 14627, USA}
\email{cmlr@math.rochester.edu}
\author{Pascal Vogt}
\address{ifb AG, Neumarkt-Galerie, D-50667 K\"{o}ln, Germany}
\email{Pascal.Vogt@ifbAG.com}

\begin{abstract}
In a recent work, Fleischmann and Mueller (2004) showed the existence of a
super-Brownian motion in $\mathbb{R}^{d}$, $d=2,3$, with extra birth at the
origin. Their construction made use of an analytical approach based on the
fundamental solution of the heat equation with a one point potential worked
out by Albeverio et al.\ (1995). The present note addresses two properties of
this measure-valued process in the three-dimensional case, namely the scaling
of the process and the large scale behavior of its mean.\vspace{-15pt}

\end{abstract}
\maketitle

\thispagestyle{empty}

\section{Introduction}

A super-Brownian motion in $\mathbb{R}$ with a single point source $\delta
_{0}$ was constructed in Engl\"{a}nder \& Fleischmann
\cite{EnglaenderFleischmann2000.prod}. It was shown that its expected mass
grows exponentially in time, and is in the mass-rescaled limit distributed in
space as $x\mapsto\mathrm{e}^{-|x|}.$ In Engl\"{a}nder \& Turaev
\cite{EnglanderTuraev2002} it is even proved that the random measures
themselves grow in law exponentially as time increases, and are otherwise in
the mass-rescaled limit spatially situated with the same shape except an
overall random factor. The probabilistic effect behind the non-trivial
existence of the model is the fact that a Brownian particle in $\mathbb{R}$
hits the origin with certainty and that it has there a non-degenerate local
time, serving as an additional birth rate for the random creation of mass.

In higher dimensions, a Brownian particle fails to hit the origin, and a local
time would degenerate. Nevertheless, Fleischmann \& Mueller
\cite{FleischmannMueller2004.birth} succeeded in constructing a super-Brownian
motion in $\mathbb{R}^{d},$ $d=2,3,$ with a single point source. They heavily
used well-known analytical facts from mathematical physics concerning Laplace
operators with one-point-potentials. Heuristically, some additional rescaling
enters the regularization of the delta function (serving as single point
source). Properties of this new super-Brownian motion are not known so far.
The purpose of the present note is to get some progress by studying its
scaling and the large scale behavior of its expectation in the
three-dimensional case.

\subsection{The heat equation with one-point-potential}

The Schr\"{o}\-dinger equation with a one-point-poten\-tial is studied in
quantum theory to describe singular electromagnetic effects on quantum
particles, see e.g.\ the monograph Albeverio et al.\ \cite[Part~I]%
{AlbeverioGesztesyHoeg-KrohnHolden1988}. By analytic continuation, solutions
to the Schr\"{o}dinger equation can be (at least formally) obtained via
solutions of the heat equation.

Formally, the heat equation with a one-point-potential is given by
\begin{equation}
\partial_{t}u=\Delta u+\delta_{0}^{(\alpha)}u=:\Delta^{\!(\alpha)}u,
\end{equation}
where $\partial_{t}$ denotes the derivative with respect to time, $\Delta$ is
the $d$--dimensional Laplacian, and $u:(0,\infty)\times\dot{\mathbb{R}}%
^{d}\rightarrow\mathbb{R}_{+}$ is a time-space field, where $\dot{\mathbb{R}%
}^{d}:=\mathbb{R}^{d}\setminus\{0\}$ with the Euclidean metric is locally
compact. If we denote by $B_{\varepsilon}(y)$ an open ball around
$y\in\mathbb{R}^{d}$ of radius $\varepsilon>0$, then having in mind that
$\,\varepsilon^{-d}\mathbf{1}_{B_{\varepsilon}(0)}\approx\delta_{0\,},\,$ the
operator $\Delta^{\!(\alpha)}:=\Delta+\delta_{0\,}^{(\alpha)}$ is
heuristically the limit as $\varepsilon\downarrow0$ of the operator%
\begin{equation}
\Delta_{\varepsilon}^{\!(\alpha)}:=\Delta+h(d,\alpha,\varepsilon
)\,\varepsilon^{-d}\mathbf{1}_{B_{\varepsilon}(0)},
\end{equation}
where $h(d,\alpha,\varepsilon)$ is some additional rescaling factor which
depends on a parameter $\alpha$ at least. Restricting to $d=3,$ the function
$h$ can be chosen as
\begin{equation}
h(3,\alpha,\varepsilon):=\tfrac{\pi^{2}}{4}\,\varepsilon-8\pi^{2}%
\alpha\varepsilon^{2},\qquad\alpha\in\mathbb{R},\ \,\varepsilon>0, \label{h}%
\end{equation}
(cf.\ \cite[(H.74)]{AlbeverioGesztesyHoeg-KrohnHolden1988}).

Physically, $\alpha$ in the case $\alpha<0$ is related to the \emph{scattering
length} $\,\mathrm{sl}_{\alpha}:=-(4\pi\alpha)^{-1}\,$ of the free Laplace
operator $\Delta$ with respect to the interaction Laplacian $\Delta
^{\!(\alpha)}$. Roughly speaking, the scattering length describes the average
distance a free particle manages to go before any interaction takes place. So,
if $\alpha\downarrow-\infty$ the scattering length $\,\mathrm{sl}_{\alpha
}\downarrow0\,$ becomes smaller and we expect more interaction. For
$\alpha\geq0$ there is \emph{no} proper physical interpretation of
$\,\mathrm{sl}_{\alpha\,}$ as the point spectrum of $\Delta^{\!(\alpha)}$ is
empty (see \cite[Theorem~I.1.4]{AlbeverioGesztesyHoeg-KrohnHolden1988}).

The \emph{fundamental solution}\/ $p^{\alpha}$ to the equation
\begin{equation}
\partial_{t}u=\Delta^{\!(\alpha)}u\quad\text{ on }(0,\infty)\times
\dot{\mathbb{R}}^{d},\quad d=2,3,
\end{equation}
which provides the basis for the analytical construction of the superprocess
in \cite{FleischmannMueller2004.birth}, have been computed in Albeverio et
al.\ \cite{AlbeverioBrzezniakDabrowski1995}. In $d=3$ (the two-dimensional
case is analytically more delicate, which is the reason we restrict to $d=3),$
the \emph{one-point-interaction heat kernel}\/ $p^{\alpha}$ is given by
\begin{equation}
p_{t}^{\alpha}(x,y)\,=\,p_{t}(x,y)+\frac{2t}{|x||y|}\,p_{t}%
\big(|x|+|y|\big)-\frac{8\pi\alpha t}{|x||y|}\int_{0}^{\infty}%
du\,\mathop {\mathrm{e}^{ -4\pi \alpha u}}p_{t}\big(u+|x|+|y|\big), \label{pa}%
\end{equation}
$t>0,$ $\,x,y\neq0,$ where $p$ is the usual \emph{free}\/ heat kernel defined
by,
\begin{equation}
p_{t}(x,y):=(4\pi t)^{-d/2}\exp\big(-|y-x|^{2}/4t\big),
\end{equation}
and with a slight abuse of notation,
\begin{equation}
p_{t}(r):=(4\pi t)^{-d/2}\exp(-r^{2}/4t),\quad\quad t>0,\quad r\geq0.
\end{equation}
Also recall the scaling of the free heat kernel, i.e. for all $k,t>0$ and
$x,y\in\mathbb{R}^{d}$,
\begin{equation}
p_{t}(x,y)=k^{d/2}p_{kt}(k^{1/2}x,k^{1/2}y). \label{scalingP}%
\end{equation}

Note, that the last term in (\ref{pa}) is always finite and disappears for
$\alpha=0$. Moreover, $\alpha\mapsto p^{\alpha}$ is pointwise continuous and
decreasing, and we have the (pointwise) convergences $p^{\alpha}%
\uparrow+\infty$ as $\alpha\downarrow-\infty$ (i.e. the fundamental solution
explodes which can be interpreted as \emph{immediate interaction}), whereas
$p^{\alpha}\downarrow p$ as $\alpha\uparrow+\infty$ leads the free case
(i.e.\ the interaction disappears).

Rigorously, the family $\{\Delta^{\!(\alpha)}:\,\alpha\in\mathbb{R}\}$ of
operators are defined as \emph{all}\/ self-adjoint extensions on the Hilbert
space $\mathcal{L}^{2}(\dot{\mathbb{R}}^{d},dx)$ of the Laplacian $\Delta$
acting on $\mathcal{C}_{\mathrm{com}}^{\infty}(\dot{\mathbb{R}}^{d}),$ the
space of unboundedly differentiable functions on $\dot{\mathbb{R}}%
^{d}=\mathbb{R}^{d}\setminus\{0\}$ with compact support (see e.g.
\cite[Chapters I.1 and I.5]{AlbeverioGesztesyHoeg-KrohnHolden1988}). Hence,
although the $p^{\alpha}$ differ from the free heat kernel $p$, they solve the
heat equation%
\begin{equation}
\partial_{t}p_{t}^{\alpha}(x,y)=\Delta p_{t}^{\alpha}(x,y)\quad\text{on
}(0,\infty)\times\dot{\mathbb{R}}^{d},
\end{equation}
with the Laplacian $\Delta$ acting either on the variable $x$ or $y.$ In
particular, $(t,x,y)\mapsto p_{t}^{\alpha}(x,y)$ is jointly continuous on
$(0,\infty)\times\dot{\mathbb{R}}^{d}\times\dot{\mathbb{R}}^{d}$%
.\thinspace\ Let us denote by $S^{\alpha}$ the semigroup associated with the
kernel $p^{\alpha}$, i.e.
\begin{equation}
S_{t}^{\alpha}\varphi(x):=\int_{\dot{\mathbb{R}}^{d}}dy\,\varphi
(y)\,p_{t}^{\alpha}(x,y). \label{sg}%
\end{equation}
Note that $S^{\alpha}$ is \emph{not} a contraction semigroup and so there is
\emph{no} stochastic process generated by this flow. The following Lemma shows
that the kernel $p^{\alpha}$ has a similar scaling behavior as the free heat
kernel $p.$

\begin{lemma}
[\textbf{Scaling of the $p^{\alpha}$}]\label{Pscaling}We have, for all
$\,k,t>0,$ $\,x,y\in\dot{\mathbb{R}}^{3},$ and $\,\alpha\in\mathbb{R}$,
\begin{equation}
p_{t}^{\alpha}(x,y)=k^{3/2}\,p_{kt}^{k^{-1/2}\alpha}(k^{1/2}x,k^{1/2}y).
\end{equation}

\end{lemma}

\begin{proof}
That follows immediately from the definition (\ref{pa}) of the $p^{\alpha}$
and the scaling (\ref{scalingP}) of the free heat kernel $p$.
\end{proof}

\subsection{The flow associated with the one-point-interaction heat kernel}

This section is devoted to introduce a space of functions $\Phi$ on which the
flow $S^{\alpha}$ acts as a strongly continuous linear semigroup (see
\cite[Section~2]{FleischmannMueller2004.birth} for details). Let $\phi$ denote
the \emph{weight and reference function}%
\begin{equation}
\phi(x):=|x|^{-1},\qquad x\in\dot{\mathbb{R}}^{3}=\mathbb{R}^{3}%
\setminus\{0\}. \label{ref.function}%
\end{equation}
For fixed $\varrho\in(1,2),$ let $\mathcal{H}=\mathcal{H}^{\varrho}$ denote
the space of measurable functions $\varphi$ on $\dot{\mathbb{R}}^{3}$ for
which
\begin{equation}
\Vert\varphi\Vert_{\mathcal{H}}:=\Big(\int_{\dot{\mathbb{R}}^{3}}%
dx\,\phi(x)\,|\varphi(x)|^{\varrho}\Big)^{1/\varrho}<\infty.
\end{equation}
Then $\big(\mathcal{H},\Vert\cdot\Vert_{\mathcal{H}}\big)$ is a Banach space,
where as usual we do not distinguish between equivalence classes and their
representatives. Now, let $\Phi=\Phi^{\varrho}$ denote the set of all
\emph{continuous}\/ functions $\varphi:\dot{\mathbb{R}}^{3}\rightarrow
\mathbb{R}$ such that $\varphi\in\mathcal{H}$ and
\begin{equation}
0\leq\varphi\leq C\,\phi\quad\text{for some constant }\,C=C_{\varphi}>0.
\end{equation}
We endow $\Phi$ with the topology inherited from $\mathcal{H}$. Note that the
set $\mathcal{C}_{\text{com}}^{+}=\mathcal{C}_{\text{com}}^{+}(\dot
{\mathbb{R}}^{3})$ of all non-negative, continuous functions on $\dot
{\mathbb{R}}^{3}$ with compact support is contained in $\Phi$. We remark that
$\varphi\in\Phi$ might have a singularity at $x=0$ of order $|x|^{-\xi}$ with
$0<\xi<1\,.$ The linear semigroup $S^{\alpha}$ introduced in (\ref{sg}) is
strongly continuous on the cone $\Phi=\Phi^{\varrho},$ cf.\ Corollary~2.12 in
\cite{FleischmannMueller2004.birth}.

\subsection{Super-Brownian motion with a single point source}

Denote by $\mathcal{M}=\mathcal{M}(\dot{\mathbb{R}}^{3})$ the set of all
(Radon) measures $\mu$ on $\dot{\mathbb{R}}^{3}$ such that $\langle\mu
,\varphi\rangle<\infty$ for all $\varphi\in\Phi$. Recalling that
$\mathcal{C}_{\text{com}}^{+}\subset\Phi,$ endow $\mathcal{M}$ with the vague topology.

Fix a constant $\eta>0$ (branching rate). Suppose $0<\beta<1$ (the finite
variance branching case $\beta=1$ has been excluded in
\cite{FleischmannMueller2004.birth} for $d=3$ for technical reasons). Then for
each $\alpha\in\mathbb{R}$, there is a non-degenerate $\mathcal{M}$--valued
(time-homogeneous) Markov process $X^{\alpha}$ such that for (deterministic)
starting measures $\mu\in\mathcal{M}\ $and for $\varphi\in\Phi$,
\begin{equation}
-\log\mathbb{E}_{\mu}\exp\langle X_{t}^{\alpha},-\varphi\rangle=\big\langle\mu
,v(t,\cdot)\big\rangle,\qquad t>0, \label{log-Lap}%
\end{equation}
where $\big\{v(t,x):\,t\geq0,\,x\in\dot{\mathbb{R}}^{3}\big\}$ is the unique
non-negative solution of the integral equation related to the $\Phi$-valued
evolution equation
\begin{equation}
\left\{
\begin{array}
[c]{c}%
\partial_{t}v=\Delta^{\!(\alpha)}v-\eta\,v^{1+\beta}\quad\text{ on }%
(0,\infty),\\[4pt]%
v(0+,\,\cdot\,)=\varphi
\end{array}
\right.
\end{equation}
(see \cite[Theorem~4.4]{FleischmannMueller2004.birth}). That is,%
\begin{equation}
v(t,x)\ =\ \int_{\dot{\mathbb{R}}^{3}}\!\mathrm{d}y\ p_{t}^{\alpha
}(x,y)\,\varphi(y)-\eta\int_{0}^{t}\!\mathrm{d}s\,\int_{\dot{\mathbb{R}}^{3}%
}\!\mathrm{d}y\ p_{t-s}^{\alpha}(x,y)\,v^{1+\beta}(s,y), \label{log-Lap.equ}%
\end{equation}
$t>0,$\thinspace\ $x\in\dot{\mathbb{R}}^{3}.$\thinspace\ Clearly, the first
moments of $X^{\alpha}$ are determined by the $S^{\alpha}$ flow to be
\begin{equation}
\mathbb{E}_{\mu}\langle X_{t}^{\alpha},\varphi\rangle=\langle\mu,S_{t}%
^{\alpha}\varphi\rangle, \label{ew}%
\end{equation}
for all starting measures $\mu\in\mathcal{M}$, $t\geq0,$ and $\varphi\in\Phi$.

\section{Large scale behavior}

\subsection{A scaling property}

\begin{proposition}
[\textbf{A scaling property}]\label{P.scaling}Let $\,t,k>0,$\thinspace
\ $\mu\in\mathcal{M},$\thinspace\ and\thinspace\ $\alpha,\lambda_{k}%
\in\mathbb{R}.$\thinspace\ Then%
\begin{equation}
\left\{  k^{-1/\beta}X_{kt}^{\lambda_{k}\alpha}(k^{1/2}\,\cdot
\,)\,\big|\,X_{0}^{\lambda_{k}\alpha}=k^{1/\beta}\mu(k^{-1/2}\,\cdot
\,)\right\}  \ \overset{\mathcal{L}}{=}\ \left\{  X_{t}^{k^{1/2}\lambda
_{k}\alpha}\,\big|\,X_{0}^{k^{1/2}\lambda_{k}\alpha}=\mu\right\}  .
\label{claim}%
\end{equation}

\end{proposition}

Of course, the cases $\,\lambda_{k}=k^{-1/2}$\thinspace\ or even $\,\alpha
=0$\thinspace\ are particularly nice.%

\proof
For $\,\varphi\in\Phi$\thinspace\ fixed,%
\begin{equation}%
\big\langle
k^{-1/\beta}X_{kt}^{\lambda_{k}\alpha}(k^{1/2}\,dy),\,\varphi%
\big\rangle
\ =\
\big\langle
X_{kt}^{\lambda_{k}\alpha},\,k^{-1/\beta}\varphi(k^{-1/2}\,\cdot\,)%
\big\rangle
,
\end{equation}
hence, by (\ref{log-Lap}) and (\ref{log-Lap.equ}),
\begin{gather}
-\log\mathbb{E}_{k^{1/\beta}\mu(k^{-1/2}\,\cdot\,)}\exp%
\big\langle
k^{-1/\beta}X_{kt}^{\lambda_{k}\alpha}(k^{1/2}\,dy),\,\varphi%
\big\rangle
\ =\ -\log\mathbb{E}_{k^{1/\beta}\mu(k^{-1/2}\,\cdot\,)}\exp%
\big\langle
X_{kt}^{\lambda_{k}\alpha},\,k^{-1/\beta}\varphi(k^{-1/2}\,\cdot\,)%
\big\rangle
\nonumber\\
=\
\big\langle
k^{1/\beta}\mu(k^{-1/2}\,\cdot\,),v(kt,\,\cdot\,)%
\big\rangle
\ =\
\big\langle
\mu,k^{1/\beta}v(kt,k^{1/2}\,\cdot\,)%
\big\rangle
\end{gather}
where $\big\{v(t^{\prime},x^{\prime}):\,t^{\prime}\geq0,\,x^{\prime}\in
\dot{\mathbb{R}}^{3}\big\}$ is the non-negative solution of the integral
equation related to the function-valued evolution equation
\begin{equation}
\left\{
\begin{array}
[c]{c}%
\partial_{t}v=\Delta^{\!(\lambda_{k}\alpha)}v-\eta\,v^{1+\beta}\quad\text{ on
}(0,\infty),\\[4pt]%
v(0+,\,\cdot\,)=k^{-1/\beta}\varphi(k^{-1/2}\,\cdot\,).
\end{array}
\right.
\end{equation}
More precisely,%
\begin{align}
k^{1/\beta}v(kt,k^{1/2}x)\ =\  &  k^{1/\beta}\int_{\dot{\mathbb{R}}^{3}%
}\!\mathrm{d}y\ p_{kt}^{\lambda_{k}\alpha}(k^{1/2}x,y)\,k^{-1/\beta}%
\varphi(k^{-1/2}y)\,\\
&  \ -\,k^{1/\beta}\eta\int_{0}^{kt}\!\mathrm{d}s\,\int_{\dot{\mathbb{R}}^{3}%
}\!\mathrm{d}y\ p_{kt-s}^{\lambda_{k}\alpha}(k^{1/2}x,y)\,v^{1+\beta
}(s,y).\nonumber
\end{align}
By a change of variable,%
\begin{align}
k^{1/\beta}v(kt,k^{1/2}x)\ =\  &  \int_{\dot{\mathbb{R}}^{3}}\!\mathrm{d}%
y\ k^{3/2}p_{kt}^{\lambda_{k}\alpha}(k^{1/2}x,k^{1/2}y)\,\varphi(y)\,\\
&  \ -\,k^{1/\beta}\eta\int_{0}^{t}\!\mathrm{d}s\,k\int_{\dot{\mathbb{R}}^{3}%
}\!\mathrm{d}y\ k^{3/2}p_{kt-ks}^{\lambda_{k}\alpha}(k^{1/2}x,k^{1/2}%
y)\,v^{1+\beta}(ks,k^{1/2}y).\nonumber
\end{align}
Hence, by Lemma~\ref{Pscaling},%
\[
k^{1/\beta}v(kt,k^{1/2}x)\ =\ \int_{\dot{\mathbb{R}}^{3}}\!\mathrm{d}%
y\ p_{t}^{k^{1/2}\lambda_{k}\alpha}(x,y)\,\varphi(y)\,-\,k^{1/\beta}\eta
\int_{0}^{t}\!\mathrm{d}s\,k\int_{\dot{\mathbb{R}}^{3}}\!\mathrm{d}%
y\ p_{t-s}^{k^{1/2}\lambda_{k}\alpha}(x,y)\,v^{1+\beta}(ks,k^{1/2}y).
\]
Since $\,1/\beta+1-(1/\beta)(1+\beta)=0$\thinspace\ we see that $\,k^{1/\beta
}v(kt,k^{1/2}x)=:w_{k}(t,x)$\thinspace\ satisfies the equation
\begin{equation}
w_{k}(t^{\prime},x^{\prime})\ =\ \int_{\dot{\mathbb{R}}^{3}}\!\mathrm{d}%
y\ p_{t^{\prime}}^{k^{1/2}\lambda_{k}\alpha}(x^{\prime},y)\,\varphi
(y)-\eta\int_{0}^{t^{\prime}}\!\mathrm{d}s\,\int_{\dot{\mathbb{R}}^{3}%
}\!\mathrm{d}y\ p_{t^{\prime}-s}^{k^{1/2}\lambda_{k}\alpha}(x^{\prime
},y)\,w_{k}^{1+\beta}(s,y),
\end{equation}
$t^{\prime}>0,$\thinspace\ $x^{\prime}\in\dot{\mathbb{R}}^{3}.$\thinspace\ By
uniqueness of solutions of the log-Laplace equation (\ref{log-Lap.equ}) and by
(\ref{log-Lap}), claim (\ref{claim}) follows.%
\endproof

\subsection{Expectation of the scaled $X^{\alpha}$}

Before we can state the result, we have to introduce some notation. The
limiting measure will be expressed by means of the kernel
\begin{equation}
\vartheta_{t}^{\alpha}(x,y):=\frac{2t}{|x|\,|y|}\,p_{t}\big(|y|\big)-\frac
{8\pi\alpha t}{|x|\,|y|}\int_{0}^{\infty}%
du\,\mathop {\mathrm{e}^{ -4\pi \alpha u}}p_{t}\big(u+|y|\big),
\label{not.vartheta}%
\end{equation}
for $\alpha\in\mathbb{R},$ $\,t>0,$ and $x,y\in\dot{\mathbb{R}}^{3}$. Note
that the integral is always finite, hence for $\alpha=0$ the second term
disappears. Moreover, the kernel $\vartheta^{\alpha}$ is always non-negative.
This holds trivially whenever $\alpha<0$, and to see this for $\alpha>0,$ use
the estimate
\begin{equation}
p_{t}\big(u+|y|\big)\leq p_{t}\big(|y|\big).
\end{equation}
We extend the definition of $\vartheta^{\alpha}$ by setting
\begin{equation}
\vartheta_{t}^{\alpha}(x,y):\equiv\left\{
\begin{array}
[c]{lll}%
0, & \text{ if } & \alpha=+\infty,\\
+\infty, & \text{ if } & \alpha=-\infty.
\end{array}
\right.
\end{equation}
The so defined kernels $\vartheta^{\alpha}$ turn out to be pointwise
continuous in $\alpha\in\lbrack-\infty,+\infty]$ (which follows from the
arguments of the proof of Theorem \ref{totalmass} below).

\begin{theo}
[\textbf{Large scale behavior of the mean}]\label{totalmass}For $t>0$%
,\thinspace\ $\alpha,\lambda_{k}\in\mathbb{R},$ and all starting measures
$X_{0}^{\alpha}=\mu\in\mathcal{M}$ satisfying $\left\langle \mu,\phi
\right\rangle <\infty$, we have the convergence in $\mathcal{M}$,
\begin{equation}
\lim_{k\uparrow\infty}k^{-1/2}\,{}\mathbb{E}_{\mu}\big[X_{kt}^{\lambda
_{k}\alpha}(k^{1/2}\,dy)\big]\,=\
\big\langle
\mu,\vartheta_{t}^{\alpha^{\ast}}(\,\cdot\,,y)%
\big\rangle
\,dy, \label{asymp}%
\end{equation}
provided that $\,\alpha^{\ast}:=\lim_{k\uparrow\infty}k^{1/2}\lambda_{k}%
\alpha\in\lbrack-\infty,+\infty]$.
\end{theo}

\begin{proof}
Fix $\varphi\in\mathcal{C}_{\mathrm{com}}^{+}(\dot{\mathbb{R}}^{3}).$ Using
formula (\ref{ew}) for the first moment of $X^{\alpha}$ and substitution, we
obtain
\begin{gather}
k^{-1/2}\,{}\mathbb{E}_{\mu}\big\langle X_{kt}^{\lambda_{k}\alpha}%
,\varphi(k^{-1/2}\,\cdot\,)\big\rangle=k^{-1/2}\int_{\dot{\mathbb{R}}^{3}}%
\mu(dx)\int_{\dot{\mathbb{R}}^{3}}dy\,p_{kt}^{\lambda_{k}\alpha}%
(x,y)\,\varphi(k^{-1/2}y)\nonumber\\
=k^{-1/2}\int_{\dot{\mathbb{R}}^{3}}\mu(dx)\int_{\dot{\mathbb{R}}^{3}%
}dy\,k^{3/2}\,p_{kt}^{\lambda_{k}\alpha}(x,k^{1/2}y)\,\varphi(y).
\end{gather}
By Lemma \ref{Pscaling} this is equal to
\begin{equation}
k^{-1/2}\int_{\dot{\mathbb{R}}^{3}}\mu(dx)\int_{\dot{\mathbb{R}}^{3}}%
dy\,p_{t}^{k^{1/2}\lambda_{k}\alpha}(k^{-1/2}x,y)\,\varphi(y). \label{theo1-5}%
\end{equation}
Inserting according to definition (\ref{pa}) of $p^{\alpha},$ we get three
terms, we will deal with separately.\medskip

$1^{\circ}$ (\emph{First term}). The first term equals,
\begin{equation}
k^{-1/2}\int_{\dot{\mathbb{R}}^{3}}\mu(dx)\int_{\dot{\mathbb{R}}^{3}}%
dy\,p_{t}(k^{-1/2}x,y)\,\varphi(y). \label{first.term}%
\end{equation}
This double integral is finite and vanishes as $k\uparrow\infty$. To see this,
let us restrict the outer integral first to $|x|>K$ where we specify $K\geq1$
later. We call this restricted integral $\,I_{K\,}.$\thinspace\ We use
$\,\varphi\leq C\phi$ (since $\varphi\in\Phi)$ and, with $S$ denoting the free
heat flow,
\begin{equation}
S_{t}\phi\leq C\phi,\qquad t\geq0, \label{heat.estimate}%
\end{equation}
with changed constant $C$ (see \cite[Lemma~2.1]{FleischmannMueller2004.birth})
to arrive at
\begin{equation}
I_{K}\ \leq\ Ck^{-1/2}\int_{|k^{-1/2}x|>K}\mu(dx)\,\phi(k^{-1/2}%
x)\ =\ C\int_{|x|>k^{1/2}K}\mu(dx)\,\phi(x) \label{first}%
\end{equation}
which can be made arbitrarily small uniformly in $k$ by choosing $K$
sufficiently large (by our assumption on $\mu)$. It remains to deal with the
case $|x|\leq K$ for fixed $K.$ We split the internal integral in
(\ref{first.term}) as follows. First, if $|k^{-1/2}x-y|\geq|y|/2$, then
\begin{equation}
p_{t}(k^{-1/2}x,y)\leq p_{t}\big(|y|/2\big),
\end{equation}
which leads to the bound
\begin{equation}
k^{-1/2}\int_{|x|\leq K}\mu(dx)\int_{\dot{\mathbb{R}}^{3}}dy\,p_{t}%
\big(|y|/2\big)\,\varphi(y)\longrightarrow0\quad\text{as }k\uparrow\infty,
\end{equation}
the $\mu(dx)$-integral is finite as $\left\langle \mu,\phi\right\rangle
<\infty$. On the other hand, if $|k^{-1/2}x-y|<|y|/2$, then $-k^{-1/2}%
|x|+|y|<|y|/2$ which implies $|y|<2k^{-1/2}|x|\leq2K$. Hence as $p_{t}%
(k^{-1/2}x,y)\leq Ct^{-3/2}$ and $\varphi\leq C\phi$, we get the upper
estimate
\begin{equation}
C_{t}\,k^{-1/2}\int_{|x|\leq K}\mu(dx)\int_{|y|<\,2K}dy\,\phi
(y)\longrightarrow0\quad\text{as }k\uparrow\infty
\end{equation}
(the $dy$-integral is finite, since we are in dimension three).\medskip

$2^{\circ}$ (\emph{Second term}). The second term reads
\begin{equation}
\int_{\dot{\mathbb{R}}^{3}}\mu(dx)\int_{\dot{\mathbb{R}}^{3}}dy\ \frac
{2t}{|x|\,|y|}\,p_{t}\big(k^{-1/2}|x|+|y|\big)\,\varphi(y)\ =:\ I\!I_{k\,}.
\label{theo1-1}%
\end{equation}
Observe that,
\begin{equation}
p_{t}\big(k^{-1/2}|x|+|y|\big)\uparrow p_{t}\big(|y|\big)\quad\text{as
}\,k\uparrow\infty.
\end{equation}
We can apply the monotone convergence theorem to obtain the limit,
\begin{equation}
\lim_{k\uparrow\infty}\ I\!I_{k\,}=\ \int_{\dot{\mathbb{R}}^{3}}\mu
(dx)\,\frac{2t}{|x|}\int_{\dot{\mathbb{R}}^{3}}dy\ \frac{1}{|y|}%
\,p_{t}(y)\,\varphi(y), \label{theo1-gw}%
\end{equation}
where finiteness follows from $\left\langle \mu,\phi\right\rangle <\infty$ and
$\varphi\leq C\phi$. Hence, this summand gives the first part of the kernel
$\vartheta^{\alpha^{\ast}}$. \medskip

$3^{\circ}$ (\emph{Third term}). It remains to insert the scaled third term
from (\ref{pa}) into (\ref{theo1-5}) which reads as
\begin{equation}
k^{-1/2}\int_{\dot{\mathbb{R}}^{3}}\mu(dx)\int_{\dot{\mathbb{R}}^{3}}%
dy\ \frac{-8\pi tk^{1/2}\lambda_{k}\alpha}{k^{-1/2}|x|\,|y|}\int_{0}^{\infty
}du\,\mathrm{e}^{-4\pi k^{1/2}\lambda_{k}\alpha u}\,p_{t}\big(u+k^{-1/2}%
|x|+|y|\big)\,\varphi(y). \label{theo1-99}%
\end{equation}
We distinguish several cases: If $\lambda_{k}\alpha=0$ for all sufficiently
large $k$, then the third term disappears and we are done. From now on assume
$\lambda_{k}\alpha\neq0$ for all $k$. Substituting $u\mapsto4\pi
k^{1/2}|\lambda_{k}\alpha|u$ into (\ref{theo1-99}) yields,
\begin{equation}
\int_{\dot{\mathbb{R}}^{3}}\mu(dx)\int_{\dot{\mathbb{R}}^{3}}dy\,\frac
{-2t\,\text{sign}(\lambda_{k}\alpha)}{|x|\,|y|}\int_{0}^{\infty}%
du\,\mathop {\mathrm{e}^{ - \text{sign}(\rho_k \alpha)\, u}}p_{t}\Big(\frac
{u}{4\pi k^{1/2}|\lambda_{k}\alpha|}+k^{-1/2}|x|+|y|\Big)\,\varphi(y).
\label{theo1-2}%
\end{equation}

Now let $k^{1/2}|\lambda_{k}\alpha|\rightarrow\infty$. We may consider a
monotone subsequence of $k^{1/2}\lambda_{k}\alpha.$ Clearly,
\begin{equation}
p_{t}\Big(\frac{u}{4\pi k^{1/2}|\lambda_{k}\alpha|}+k^{-1/2}%
|x|+|y|\Big)\uparrow p_{t}\big(|y|\big)\quad\text{as }k\uparrow\infty,
\end{equation}
and by monotone convergence the expression (\ref{theo1-2}) converges along the
subsequence to
\begin{equation}
-\int_{\dot{\mathbb{R}}^{3}}\mu(dx)\,\frac{2t\,\text{sign}(\alpha^{\ast}%
)}{|x|}\int_{0}^{\infty}du\,\mathrm{e}^{-\text{sign}(\alpha^{\ast})u}%
\int_{\dot{\mathbb{R}}^{3}}dy\ \frac{1}{|y|}\,p_{t}(y)\,\varphi(y),
\end{equation}
which is independent of the choice of the subsequence. Note that
\begin{equation}
\text{sign}(\alpha^{\ast})\int_{0}^{\infty}du\,\mathrm{e}^{-\text{sign}%
(\alpha^{\ast})\,u}=\left\{
\begin{array}
[c]{cl}%
1 & \text{ if }\ \text{sign}(\alpha^{\ast})=1,\vspace{4pt}\\
+\infty & \text{ if }\ \text{sign}(\alpha^{\ast})=-1.
\end{array}
\right.
\end{equation}
In the first case the second and the third limiting terms cancel.

Next, let $k^{1/2}\lambda_{k}\alpha\rightarrow0$. Note, that
\begin{equation}
p_{t}\Big(\frac{u}{4\pi k^{1/2}|\lambda_{k}\alpha|}+k^{-1/2}|x|+|y|\Big)\leq
p_{t}\Big(\frac{u}{4\pi k^{1/2}|\lambda_{k}\alpha|}%
\Big)\mathop {\mathrm{e}^{ -|y|^2/4t}}. \label{split}%
\end{equation}
In this case the double integral (\ref{theo1-2}) is in absolute value bounded
by
\begin{equation}
\int_{\dot{\mathbb{R}}^{3}}\mu(dx)\,\int_{\dot{\mathbb{R}}^{3}}dy\,\frac
{2t}{|x|\,|y|}\mathop {\mathrm{e}^{ -|y|^2/4t}}\varphi(y)\int_{0}^{\infty
}du\,\mathop {\mathrm{e}^{ u}}p_{t}\Big(\frac{u}{4\pi\,k^{1/2}|\lambda
_{k}\alpha|}\Big),
\end{equation}
which tends to 0 as $k\uparrow\infty$ as the $\mu(dx)$ and $dy$-integrals are
finite and in the $du$-integral the $p_{t}$-term compensates the
$\mathop {\mathrm{e}^{ u}}$.

It remains to deal with the case $k^{1/2}\lambda_{k}\alpha\rightarrow
\alpha^{\ast}\in\dot{\mathbb{R}}^{1}$. Note, that we only have to justify to
change the limit and integration in (\ref{theo1-2}), as substituting
$u\mapsto\big(4\pi|\alpha^{\ast}|\big)^{-1}u$ leads to the desired expression.
To justify the interchange, we estimate as in (\ref{split}). The resulting
$\mu(dx)$ and $dy$-integrals are independent of $k$ and finite, whereas to
dominate in the second integral we use $k^{1/2}|\lambda_{k}\alpha|\leq
|\alpha^{\ast}|+1$ for all sufficiently large $k$.
\end{proof}

\begin{rem}
[\textbf{Large scale total mass}]\label{R.tot.mass}Taking $\lambda_{k}\equiv1$
and choosing formally $\varphi=1$ as test function in (\ref{asymp}) yields%
\begin{equation}
\lim_{k\uparrow\infty}k^{-1/2}\,{}\mathbb{E}_{\mu}\left\langle X_{kt}^{\alpha
},1\right\rangle \ =\ \left\{
\begin{array}
[c]{ll}%
0 & \text{if }\,\alpha>0,\\
2t\left\langle \mu,\phi\right\rangle
{\displaystyle\int_{\mathbb{R}^{3}}}
dy\,\dfrac{1}{|y|}\,p_{t}(y) & \text{if }\,\alpha=0,\\
\infty & \text{if }\,\alpha<0.
\end{array}
\right.  \label{tot.mass}%
\end{equation}
A rigorous argument can be given along the lines of the previous
proof.\hfill$\Diamond$
\end{rem}

\subsection{Discussion and open problems}

Let us comment on the \emph{three cases}\/ $\alpha^{\ast}=+\infty,$
$\alpha^{\ast}\in\mathbb{R},$ and $\alpha^{\ast}=-\infty$ in
Theorem~\ref{totalmass}. In the first case, the limiting mass disappears, more
precisely, the scaled expression $\mathbb{E}_{\mu}\big[X_{kt}^{\lambda
_{k}\alpha}(k^{1/2}\,dy)\big]$ is of order $o(k^{1/2}).$ Roughly speaking, if
$\alpha^{\ast}=+\infty,$ then there are no interactions in the scaling limit
(free case). In the second case, $\alpha^{\ast}\in\mathbb{R},$ the former
expectation is about $k^{1/2}\,%
\big\langle
\mu,\vartheta_{t}^{\alpha^{\ast}}(\,\cdot\,,y)%
\big\rangle
\,dy.$ Note that these measures are decreasing in $\alpha^{\ast}.$ Finally, if
$\alpha^{\ast}=-\infty,$ we have immediate interaction in the large scale
limit leading to the explosion of the expected mass.\smallskip

Clearly, to describe only the large scale behavior of the \emph{expected
processes}\/ is unsatisfactory. It is desirable to get something similar for
the processes themselves. Recall that in the one-dimensional case the large
time behavior of the process itself is known from \cite{EnglanderTuraev2002}.
However, we stress the fact, that the process in three dimensions is expected
to have quite different features. For instance, if $\,\alpha=0,$%
\thinspace\ then according to Remark~\ref{R.tot.mass} the total mass grows
with a power order, whereas in one dimension the growth is exponential.
Moreover, in the three-dimensional case one needs additionally to contract the
normalized measures to get a limit. For the measures themselves, scaled as in
Theorem~\ref{totalmass}, \emph{there might be extinction in law}\/ despite
convergence of their expectations.\smallskip

Another \emph{open problem}\/ is the large scale behavior of $\mathbb{E}%
X^{\alpha}$ in the two-dimensional case, in which the fundamental solutions
$p^{\alpha}$ from \cite{AlbeverioBrzezniakDabrowski1995} are analytically more
delicate, see e.g. \cite[formula~(2.30)]{FleischmannMueller2004.birth}. In
particular, a scaling property as in Lemma~\ref{Pscaling} is not
available.\bigskip

\noindent\emph{Acknowledgement}. We thank Hagen Neidhardt of the WIAS for a
helpful discussion on semigroups.

\bibliographystyle{alpha}
\bibliography{bibtex,bibtexmy}

\newcommand{\noopsort}[1]{}
\begin{thebibliography}{AGHKH88}

\bibitem[ABD95]{AlbeverioBrzezniakDabrowski1995}
S.~Albeverio, Z.~Brze{\'z}niak, and L.~Dabrowski.
\newblock Fundamental solution of the heat and {S}chr\"odinger equations with
  point interaction.
\newblock {\em J. Funct. Anal.}, 130:220--254, 1995.

\bibitem[AGHKH88]{AlbeverioGesztesyHoeg-KrohnHolden1988}
S.~Albeverio, F.~Gesztesy, R.~H{\o}egh-Krohn, and H.~Holden.
\newblock {\em Solvable Models in Quantum Mechanics}.
\newblock Springer-Verlag, New York, 1988.

\bibitem[EF00]{EnglaenderFleischmann2000.prod}
J.~Engl{\"a}nder and K.~Fleischmann.
\newblock Extinction properties of super-{B}rownian motions with additional
  spatially dependent mass production.
\newblock {\em Stoch. Proc. Appl.}, 88(1):37--58, 2000.

\bibitem[ET02]{EnglanderTuraev2002}
J.~Engl{\"a}nder and D.~Turaev.
\newblock A scaling limit theorem for a class of superdiffusions.
\newblock {\em Ann. Probab.}, 30(2):683--722, 2002.

\bibitem[FM04]{FleischmannMueller2004.birth}
K.~Fleischmann and C.~Mueller.
\newblock Super-{B}rownian motion with extra birth at one point.
\newblock {\em SIAM J. Math. Analysis}, 36(3):740--772, 2004.

\end{thebibliography}



\end{document}